\documentclass[a4paper,12pt]{article}
\usepackage{tabularx} 
\usepackage{amsmath,amsthm,verbatim,amssymb,amsfonts,amsopn,amscd}
\usepackage{graphicx} 
\usepackage[margin=1in,a4paper]{geometry} 

\usepackage[square,numbers]{natbib}
\usepackage[utf8]{inputenc}
\usepackage[T1]{fontenc}

\usepackage[UKenglish]{babel}
\usepackage{enumerate}
\usepackage{paralist} 
\usepackage[final]{hyperref} 
\hypersetup{
	colorlinks=true,       
	linkcolor=blue,        
	citecolor=blue,        
	filecolor=magenta,     
	urlcolor=blue         
}

\usepackage{tikz}

\theoremstyle{plain}
\newtheorem{theorem}{Theorem}[section]

\theoremstyle{definition}

\DeclareMathOperator{\curl}{curl}

\newcommand{\IN}{\mathbb{N}}

\newcommand{\IR}{\mathbb{R}}

\renewcommand{\vec}[1]{\overrightarrow{#1}}
\newcommand{\la}{\langle}
\newcommand{\ra}{\rangle}

\title{The Fundamental Theorem of Calculus in higher dimensions}
\author{\href{mailto:filip.bar.research@gmail.com}{Filip B\'{a}r}}
\date{22nd February 2024}

\begin{document}

\maketitle
\begin{abstract}
  \noindent We generalise the Fundamental Theorem of Calculus to higher dimensions. Our generalisation is based on the observation that the antiderivative of a function of $n$-variables is a solution of a partial differential equation of order $n$ generalising the classical case. The generalised Fundamental Theorem of Calculus then states that the $n$-dimensional integrals over $n$-dimensional axis-parallel rectangular hypercuboids is given by a combinatorial formula evaluating the antiderivative on the vertices of the hypercuboid. 
\end{abstract}

\section{Introduction}\label{sec:intro}

\noindent This note is about generalising the Fundamental Theorem of Calculus to higher dimensions. Contrary to the widespread belief that the generalisation is given by the Stokes-Cartan theorem, our approach is based on generalising the notion of an antiderivative to functions of $n$-variables yielding a very different result. Our main result states that for continuous functions of $n$-variables such an antiderivative always exists, and that we can use it to evaluate $n$-dimensional integrals over $n$-dimensional axis-parallel rectangular hypercuboids with a simple combinatorial formula.

The proof of the main result turns out to be a simple induction proof based on Fubini and the classical Fundamental Theorem of Calculus. In addition, the combinatorial formula does not generalise beyond (diffeomorphic images of) parallelotopes. Despite the basic character of the main result and these constraints, we believe that the Fundamental Theorem in $n$-dimensions is of great conceptual significance, nonetheless. Notably, it leads us to an extension of Cartan's theory of differential forms linking the integrability of certain PDEs to the geometry of the underlying manifold \cite{Bar:extending_Cartan, Bar:towards_synthetic_theory_of_integration}.

However, our primary objective here is the exposition of the generalised Fundamental Theorem of Calculus as an independent result of Calculus. We shall present the theory of differential forms elsewhere. It is exactly because it is such an elementary consequence of one of the most important theorems in Calculus that it should be more widely known. Intriguingly, our search across calculus and analysis textbooks has yielded no mention of this result thus far. Only in the context of the theorem of Schwarz and its proof via Fubini does the 2D version of the Fundamental Theorem makes an apperance as an intermediate step (see e.g. \cite{Aksoy&Martelli:partial_derivatives_Fubini}, \cite[p. 143]{Axler:measure_integration_real_analysis}). The author is therefore grateful for any pointers and references to the literature where this result has been proven, as well as for ideas or mentions of its applications in other fields.  


\section{The Fundamental Theorem of Calculus}\label{sec:FTC}

\noindent What is the generalisation of the Fundamental Theorem of Calculus? This question is typically answered by stating the Stokes-Cartan theorem for differential forms
$$
  \int_{M}d\omega = \int_{\partial M} \omega 
$$
Here $M$ denotes a smooth oriented manifold of dimension $n$ with an oriented boundary $\partial M$, and $\omega$ is a differential $n$-form. Indeed, for $n=1$ and $M=[a,b]$ this reduces to the well-known Fundamental Theorem of Calculus
$$
\int_{a}^{b}f(x)dx = F(b) - F(a) 
$$
where $F$ is the \emph{antiderivative}, or \emph{primitive} of $f$, i.e. $F'=f$. The generalisation to Stokes-Cartan proceeds as follows: The integrand $f(x)dx$ is considered as a differential $1$-form of class $C^0$, and the antiderivative $F$ as a primitive of this form; a $0$-form, which is just a real-valued function. Taking the difference on the right-hand side becomes taking the integral over the boundary of $M$, and differentiation has to be generalised to the exterior derivative. 

A major mismatch of this generalisation with the actual practice of solving integrals is revealed when we study integrals in 2D over a rectangular domain, say.
$$
  \int_{[a_1,b_1]\times[a_2,b_2]}f(x,y)\,dx\,dy
$$
Considering $f(x,y)\,dx\,dy$ as a $2$-form $f(x,y)\,dx\wedge dy$ Stokes-Cartan simplifies this integral to a line integral over the boundary of the rectangle only if $f(x,y)$ can be represented as the dot product of two vector fields, of which one is the curl of another vector field $\vec{X}$ 
$$
  f(x,y) = \la\curl \vec{X}(x,y), \vec{Y}(x,y)\ra
$$      
Most importantly, even if we can find such a representation, this is not how we solve such integrals in practice; instead, we use Fubini's theorem to represent the integral as an iterated integral
$$
  \int_{[a_1,b_1]\times[a_2,b_2]}f(x,y)\,dx\,dy = \int_{a_2}^{b_2}\int_{a_1}^{b_1}f(x,y)\,dx\,dy
$$
Applying Fubini to solve the integral requires finding an antiderivative of $f(-,y)$ first and then integrate the resulting function with respect to $y$; i.e. finding yet another antiderivative. The resulting function $F:[a_1,b_1]\times[a_2,b_2]\to \IR$ is a solution of the second-order partial differential equation 
\begin{equation}\label{eq:PDE-2D}
  \partial_x\partial_y F(x,y) = f(x,y)
\end{equation}
The function $F$ deserves to be called an \emph{antiderivative} of $f$. We can use it to evaluate the integral by taking the alternating sum of $F$ at the vertices of the rectangle:  
\begin{equation}\label{eq:FTC-2D}
  \int_{[a_1,b_1]\times[a_2,b_2]}f(x,y)\,dx\,dy = F(b_1,b_2) - F(a_1,b_2) + F(a_1,a_2) - F(a_1,b_2)
\end{equation}

The sign-rule for the differences on the right-hand side is determined by the vertices of the rectangle as follows:

  
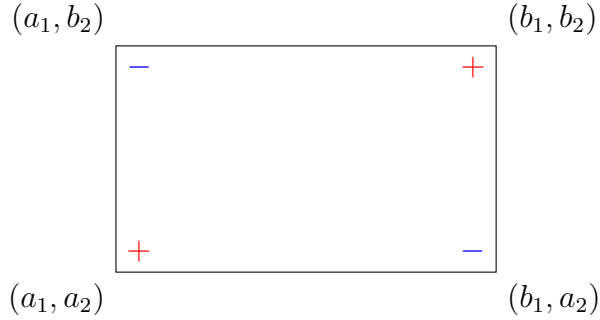
\begin{center}
    \begin{tikzpicture}
      \draw (0,0) rectangle (5,3);
      
      \node at (0,0) [below left] {$(a_1,a_2)$};
      \node at (5,0) [below right] {$(b_1,a_2)$};
      \node at (5,3) [above right] {$(b_1,b_2)$};
      \node at (0,3) [above left] {$(a_1,b_2)$};
      
      \node at (0,0) [above right, red] {$+$};
      \node at (5,0) [above left, blue] {$-$};
      \node at (5,3) [below left, red] {$+$};
      \node at (0,3) [below right, blue] {$-$};
    \end{tikzpicture}    
  \end{center}  
In line with the additivity of the integral the combinatorial formula (\ref{eq:FTC-2D}) satisfies compositionality with respect to the subdivision of the rectangle into smaller rectangles, as can be seen in the subsequent diagram.
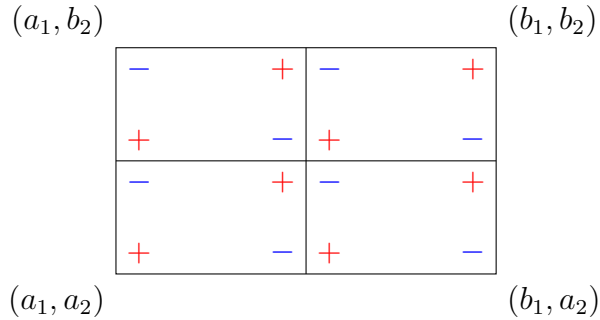
\begin{center}
  \begin{tikzpicture}
    \draw (0,0) rectangle (5,3);
    
    \draw (0,1.5) -- (5,1.5);
    \draw (2.5,0) -- (2.5,3);
    
    \node at (0,0) [below left] {$(a_1,a_2)$};
    \node at (5,0) [below right] {$(b_1,a_2)$};
    \node at (5,3) [above right] {$(b_1,b_2)$};
    \node at (0,3) [above left] {$(a_1,b_2)$};
    
    \node at (0,0) [above right, red] {$+$};
    \node at (5,0) [above left, blue] {$-$};
    \node at (5,3) [below left, red] {$+$};
    \node at (0,3) [below right, blue] {$-$};
    
    \node at (2.5,0) [above left, blue] {$-$};
    \node at (2.5,0) [above right, red] {$+$};
    
    \node at (2.5,1.5) [below right, blue] {$-$};
    \node at (2.5,1.5) [below left, red] {$+$};
    \node at (2.5,1.5) [above left, blue] {$-$};
    \node at (2.5,1.5) [above right, red] {$+$};
    
    \node at (2.5,3) [below right, blue] {$-$};
    \node at (2.5,3) [below left, red] {$+$};
    
    \node at (0,1.5) [below right, blue] {$-$};
    \node at (0,1.5) [above right, red] {$+$};
    
    \node at (5,1.5) [above left, blue] {$-$};
    \node at (5,1.5) [below left, red] {$+$};
    
  \end{tikzpicture}    
\end{center}  
Generalising the 1D case, an antiderivative of $f$ can be obtained from the iterated integral as a function of the upper bounds
$$
  F(x,y) = \int_{a_2}^{y}\int_{a_1}^{x}f(u,v)\,du\,dv
$$
However, like in the 1D case any solution of the above PDE (\ref{eq:PDE-2D}) can be used in the combinatorial formula (\ref{eq:FTC-2D}). In the 2D case the antiderivative is not unique up to a constant only, but up to a sum of two functions dependent on only one variable $x$ or $y$, respectively. This is a consequence of the second-order linear PDE (\ref{eq:PDE-2D}), and the general solution of the respective homogeneous PDE being the sum of two differentiable functions of one variable.

The formula (\ref{eq:FTC-2D}) has a direct generalisation to $n$-dimensional integrals over $n$-dimensional axis-parallel rectangular hypercuboids $I\subset \IR^n$. To be able to state it we need to introduce some notation first.

Let $2=\{0,1\}$, then $2^n$ is the set of all binary numbers of length $n$. The hypercuboid $I$ is the cartesian product of $n$ intervals $[a_j,b_j]$, with $a_j\leq b_j$, $1\leq j\leq n$. Using the pointwise order on $\IR^n$ we order the vertices of $I$ lexicographically and label them with binary numbers $b\in 2^n$ respecting the lexicographic order. Finally, let $\#_0: 2^n\to \IN$ map each binary number $b$ to the number of zeros in $b$. For example, the vertex $\vec{a}=(a_1,\ldots,a_n)$ is labelled with $0\ldots 0$, and the vertex $\vec{b}=(b_1,\ldots,b_n)$ is labelled with $1\dots 1$; hence $\#_0(0\ldots 0)=n$ and $\#_0(1\ldots 1)=0$. 

\begin{theorem}[Fundamental Theorem of Calculus] \label{thm:FTC-nD} Let $I=\prod_{j=1}^{n} [a_j,b_j]$, $\mathring{I}=\prod_{j=1}^{n} (a_j,b_j)$ and $f:I\to \IR$ a continuous function. 
\begin{enumerate}[(1)]
  \item The function 
  \begin{equation}\label{eq:antiderivative-as-integral-nD}
    F:I\to \IR, \qquad (x_1,\ldots,x_n)\mapsto \int_{a_n}^{x_n}\ldots\int_{a_1}^{x_1}f(u_1,\ldots,u_n)\,du_1\ldots du_n 
  \end{equation}
  is continuous on $I$, the mixed partial derivative $\partial_1\partial_{2}\cdots\partial_n F$ exists on $\mathring{I}$, is continuous, and a solution of the partial differential equation 
  \begin{equation}\label{eq:antiderivative-nD}
    \partial_1\cdots\partial_n F(x_1,\ldots, x_n) = f(x_1,\ldots, x_n)
  \end{equation}

  \item Let $F\in C(I,\IR)$. If the mixed partial derivative $\partial_1\cdots\partial_n F$ exists on $\mathring{I}$ and $F$ satisfies the partial differential equation (\ref{eq:antiderivative-nD}), then 
  \begin{equation}\label{eq:FTC-nD}
    \int_{I}f(x_1,\ldots,x_n)\,dx_1\ldots dx_n = \sum_{b\in 2^{n}}(-1)^{\#_0(b)}F(P_b),
  \end{equation}
  where the $P_b$ denote the lexicographically ordered vertices of $I$. Any such function $F$ will be called an \emph{antiderivative} of $f$.
\end{enumerate}
\end{theorem}

\begin{proof} The proof is a straight-forward induction on $n$. The case $n=1$ is the classical Fundamental Theorem of Calculus. For $n>1$ we can apply Fubini's theorem to represent the integral as an iterated integral. 
$$
  \int_{I}f(x_1,\ldots,x_n)\,dx_1\ldots dx_n = \int_{a_n}^{b_n}\left(\int_{I'}f(x_1,\ldots,x_n)\,dx_1\ldots dx_{n-1}\right)dx_n
$$
The inner integral is an $(n-1)$-dimensional integral over an $(n-1)$-dimensional hypercuboid $I'$ with a free parameter $x_n$, and the induction hypothesis applies. Let $G(-,x_n): I'\to \IR$ be an antiderivative of $f(-,x_n)$ defined via the iterated integral in (\ref{eq:antiderivative-as-integral-nD}), then the inner integral evaluates to
$$
  \int_{I'}f(x_1,\ldots,x_n)\,dx_1\ldots dx_{n-1} = \sum_{b\in 2^{n-1}}(-1)^{\#_0(b)}G(P_{b},x_n) 
$$
Since $f$ is continuous in $x_n$, the function $G$ is continuous in $x_n$ and thus on $I$. In particular, we can integrate each $G(x_1,\ldots,x_{n-1},-)$ over $[a_n,b_n]$:   
\begin{align*}
  \int_{I}f(x_1,\ldots,x_n)\,dx_1\ldots dx_n &= \sum_{b\in 2^{n-1}}(-1)^{\#_0(b)}\int_{a_n}^{b_n} G(P_{b},x_n)\,dx_n \\
  &= \sum_{b\in 2^{n-1}}(-1)^{\#_0(b)}(F(P_{b},b_n) - F(P_{b},a_n))  
\end{align*}

In the second step we have applied the classical Fundamental Theorem of Calculus obtaining $F(x_1,\ldots,x_{n-1},-)$ as
$$
  [a_n,b_n]\to\IR,\qquad x_n\mapsto \int_{a_n}^{x_n}G(x_1,\ldots,x_{n-1},u_n)\,du_n
$$
for each $(x_1,\ldots,x_{n-1})\in I'$. The latter also guarantees that $F$ is continuous in the $n$th argument on $[a_n,b_n]$, that $\partial_n F(x_1,\ldots,x_{n-1},-)=G(x_1,\ldots,x_{n-1},-)$ exists, and that it is continuous on $(a_n,b_n)$ for all $x_1,\ldots,x_{n-1}\in I'$. Together with $G$ being continuous on $I'$ this establishes the continuity of $F$ on $I$. Since the mixed partial derivative $\partial_{1}\cdots\partial_{n-1} G$ exists on $\mathring{I'}$, is continuous and $\partial_n F = G$, the same holds true for $\partial_1\cdots\partial_n F$ on $\mathring{I'}$. Moreover, we see that $F$ satisfies the partial differential equation (\ref{eq:antiderivative-nD}) on $\mathring{I}$ by construction and the induction hypothesis. In particular, $\partial_1\cdots\partial_n F$ is continuous on $\mathring{I}$. 

Note that each vertex of $I$ is either of the form $(P_b,b_n)$ or $(P_b,a_n)$ for a $b\in 2^{n-1}$. The vertices of the first type are labeled by binary numbers $(1,b)\in 2^{n}$, so the number of zeros in $(1,b)$ is equal to $\#_0(b)$. The vertices of the second type are labeled by binary numbers $(0,b)$, so the number of zeros in $(0,b)$ is equal to $\#_0(b)+1$. The right hand side of the last equation above can thus be written as the asserted alternating sum.

Finally, we note that any $\tilde{F}:I\to\IR$ satisfying the conditions stated in (2) will differ from $F$ by a sum of functions on $I$ that are constant in one or more of the arguments. This follows from yet another induction argument over $n$ by applying the classical Fundamental theorem of Calculus to $\partial_1\cdots\partial_n (F-\tilde{F})=0$, i.e. evaluating the iterated integral  
$$
  \int_{a_n}^{x_n}\ldots\int_{x_1}^{x_1}\partial_1\cdots\partial_n (F-\tilde{F})(u_1,\ldots,u_n)\,du_1\ldots du_n = 0
$$
for each $(x_1,\ldots,x_n)\in I$. W.l.o.g. consider a function $C:I\to \IR$ that is constant in $x_n$, and apply the combinatorial formula (\ref{eq:FTC-nD}). As argued previously, we can write (\ref{eq:FTC-nD}) as
$$
\sum_{b\in 2^{n-1}}(-1)^{\#_0(b)}(C(P_{b},b_n) - C(P_{b},a_n)),
$$
where the $P_b$ are the vertices of the $(n-1)$-dimensional hypercuboid $I'$. Since $C$ is constant in $x_n$ this sum evaluates to zero.

Similarly, if $C$ is constant in $x_j$ for some $1\leq j\leq n$, then we can partition the set of vertices of $I$ into two subsets of vertices of $(n-1)$-dimensional hypercuboids with vertices $P_b$ that have the same $x_j$-coordinate. The vertices in the first set have as the $j$th coordinate $a_j$, and the vertices in the second set $b_j$. Indeed, there are exactly two such hypercuboids and they form two opposing faces of $I$. We can write (\ref{eq:FTC-nD}) as an alteranting sum of differences over the $2^{n-1}$ vertices of both faces as above. Like in the case of $j=n$ this sum evaluates to zero for $C$ is constant in $x_j$. We conclude that both the combinatorial formula applied to $F$ and $\tilde{F}$ agree as claimed.
\end{proof}

Note that the compositionality of (\ref{eq:FTC-nD}) with respect to the division of the hypercuboid into smaller hypercuboids holds true in the $n$-dimensional case as well, which is in line with the additivity of the integral, once again. 

Next, we wish to study generalisations of (\ref{eq:FTC-nD}) to more general domains. Let $\phi: I\to U$ be a diffeomorphism\footnote{More precisely, $\phi$ is assumed a homeomorphism and $\phi|\mathring{I}$ is a diffeomorphism with image $\mathring{U}$, the interior of $U$.}. The change of variable theorem 
$$
  \int_{I}f(\phi(x_1,\ldots,x_n))
  \left|\det\partial\phi(x_1,\ldots,x_n)\right|\,dx_1\ldots dx_n = \int_{\phi(I)}f(y_1,\ldots,y_n)\,dy_1\ldots dy_n
$$
generalises (\ref{eq:FTC-nD}) directly to $n$-dimensional parallelelotopes that are not necessarily axis-parallel (and beyond). The linear transformation can be obtained as follows: Let $I$ be the unit $n$-hypercuboid and let $T$ be the change-of-basis matrix from the standard basis in $\IR^n$ to the basis of vectors spanning the parallelotope, then $\phi = T/\det T$ is the volume preserving linear transformation mapping $I$ to the parallelotope $U$. If $F$ is the antiderivative of $f\circ\varphi$, then $F\circ\phi^{-1}$ plays the role of the antiderivative of $f$ on $U$ in the sense that the integral can be evaluated by taking the alternating sum of $F\circ\phi^{-1}$ at the vertices of $U$ as in (\ref{eq:FTC-nD}). 

However, in the case of an arbitrary parallelotope it is not clear how to order the vertices anymore. This can be resolved by working with \emph{marked parellelotopes}, where the marked vertex $P_{2^n}$ equals to $\phi(1,\ldots,1)$. The signs of the other vertices $P_j$ are implied by $P_{2^n}$ as $(-1)^{d(P_j,P_{2^n})}$ where $d(P_j,P_{2^n})$ is the graph distance of vertex $P_j$ from $P_{2^n}$ on the parallelotope considered as an undirected graph. We obtain the subsequent formula for the integral over a parallelotope $U$ with marked vertex $P_{2^n}$:
\begin{equation}\label{eq:FTC-nD-parallelotope}
  \int_{U}f(x_1,\ldots,x_n)\,dx_1\ldots dx_n = \sum_{j=1}^{2^{n}}(-1)^{d(P_j,P_{2^n})}F(\phi^{-1}(P_j))
\end{equation}
The ordering of the vertices is not important, only the graph distance of the vertices from the marked vertex $P_{2^n}$ matters. By using the change-of-variable theorem we loose compositionality, though. This means the combinatorial formula (\ref{eq:FTC-nD-parallelotope}) does not extend to polyhedra that admit a tiling into parallelotopes. On the upside, (\ref{eq:FTC-nD-parallelotope}) holds also for non-linear diffeomorphisms; that is, if we parametrise a manifold with (a topological) boundary.

We should also remark that a combinatorial formula similar to (\ref{eq:FTC-nD}) cannot be expected for triangles and thus $n$-simplices for $n\geq 2$, in general. We show first how and when such a formula for the triangle can be obtained from (\ref{eq:FTC-nD-parallelotope}) for a parallelogram. We know that a triangle $PQR$ can be extended into a parallelogram $PQSR$ by reflecting it twice: first over one of its sides, say $QR$, and then a second time over the perpendicular bisector of the chosen side $QR$. Since the fourth vertex $S$ of the resulting parallelogram is an affine combination of the vertices of the triangle, namely $S = Q+R-P$, one could use the parallelogram formula to obtain a formula for the triangle taking inspiration from the method of mirror charges in electrostatics. 

This, however, would require the restriction $f|QR$ of the integrand $f$ to the chosen side to be symmetric with respect to the midpoint. Only then would the two reflections lead to a continuous extension of $f$ to the parallelogram $PQSR$. In this case the combinatorial formula would yield
$$
  \int_{PQR}f(x,y)\,dx\,dy = \frac{1}{2}(F(\phi^{-1}(Q+R-P)) - F(\phi^{-1}(Q)) + F(\phi^{-1}(P)) - F(\phi^{-1}(R))),
$$
where $F$ is the antiderivative of $f\circ\phi$ and $\phi$ is the linear isomorphism mapping the unit square $I$ to the parallelogram $PQSR$ such that $\phi(1,1)=S$. 

Note that any alternating sum over the vertices of the triangle $PQR$ will not be composable. This is because a triangle considered as an undirected graph is not bipartite and thus cannot be two-coloured. In particular, we will always have two adjacent vertices having the same sign. Since we can always divide a rectangle into two right triangles any alternating sum over the vertices of a triangle would need to imply (\ref{eq:FTC-2D}); but that is not possible. Besides introducing a factor of 2 for two of the vertices when adding the sums for each of the triangles, there is also the problem that we have two such triangulations of the rectangle. In fact, any labeling of the vertices of both triangles by plus and minus will not recover (\ref{eq:FTC-2D}).

\newpage

\section{Conclusion and Outlook}\label{sec:Conclusion}

Our exploration has revealed that contrary to the widespread belief in the Stokes-Cartan theorem as the generalisation of the Fundamental Theorem of Calculus, there exists a more direct extension rooted in the generalisation of the antiderivative concept to functions of multiple variables. Through this approach, we derived a combinatorial formula (\ref{eq:FTC-nD}) for integrals over hypercuboids, which generalises the classical formula for evaluating 1D integrals with an antiderivative.

Since the formula does not generalise well beyond parallelelotopes (respectively, the diffeomorphic images of hypercuboids), it seems to be of very limited significance when compared to the Stokes-Cartan theorem at this stage. Such a verdict would be a premature conclusion, though. In contrary, we believe that it is of great conceptual significance! 

What justifies this belief? So far we have been overly focused on the combinatorial formula (\ref{eq:FTC-nD}) itself. What we didn't do is to try to generalise the PDE defining an antiderivative to other domains than axis-parallel hypercuboids. The change-of-variable theorem gives us a first hint that this is possible. Indeed, given a diffeomorphism $\phi:I\to U$ and an antiderivative $F$ of $f\circ\phi$ as in section~\ref{sec:FTC}, the transformed PDE on $U$ has $F\circ \phi^{-1}$ as a solution, which is then used to evaluate the integral of $f$ over $U$. What if we were to integrate $f$ directly over $U$ by finding an antiderivative without refering to an apriori parametrisation $\phi$?  

Answering this question necessitates a notion of a geometric object on $U$ that can be integrated, possibly augmented with an additional geometric structure. This geometric object in question turns out to be a new type of differential form extending the classical theory of Cartan \cite{Cartan:diffforms, Cartan:diffforms_book}. This new type of differential form is best studied within the framework of Synthetic Differential Geometry \cite{Kock:Synthetic_Geometry_Manifolds, Kock:SDG}, where it is possible to model the intuition that differential $n$-forms are integrals over infinitesimal $n$-parallelelotopes directly. The combinatorial formula (\ref{eq:FTC-nD}) then becomes an example of such a differential form when applied to an infinitesimal parallelotope. In addition, it leads us to a new type of exterior derivative that, unlike Cartan's exterior derivative, doesn't square to zero. With its help it is possible for these new type of $n$-forms to have a function as a primitive for $n>1$.

The theory of these new types of differential forms is currently work in progress \cite{Bar:extending_Cartan}. A first exposition of the constructions involved in the case of two-forms can be found in \cite{Bar:towards_synthetic_theory_of_integration}.

As stated in the introduction the author is grateful for any pointers and references to the literature where the generalisation of the Fundamental Theorem has been proven, as well as for ideas or mentions of its applications in other fields.

\bibliography{References_with_links}

\end{document}